\documentclass[amstex,11pt,a4paper]{article}

\usepackage[english,russian]{babel}
\usepackage{amsfonts,amssymb}
\usepackage{graphicx}
\usepackage{latexsym,exscale}
\usepackage{amsmath}
\usepackage{amsthm}
\usepackage{blindtext}

\usepackage{amsmath}
\usepackage{relsize}

\newtheorem{Lemma}{Lemma}

\newtheorem*{myth}{Theorem}
\newtheorem{Corollary}{Corollary}

\newtheorem{Observation}{Observation}
\newtheorem{Example}{Example}

\newtheorem*{Main Theorem}{Main Theorem}

\newtheorem*{UTW}{Lemma
{\rm\cite{[KM17]}}}

\newtheorem*{ThKM}{Theorem
{\rm\cite{[KM17]}}}

\newenvironment{Proof} 
{\par\noindent{\bf{Proof.}}} 
{\hfill$\scriptstyle\blacksquare \\$}

\newenvironment{ProofMT} 
{\par\noindent{\bf{Proof of the Main Theorem.}}} 
{\hfill$\scriptstyle\blacksquare \\$}

\newcounter{Remark}

\title{О свободном разложении вербально замкнутых подгрупп в свободных произведениях конечных групп}
\author{А.М. Мажуга, А.А. Клячко}
\date{}

\makeatletter
\let\newtitle\@title
\let\newauthor\@author
\let\newdate\@date
\makeatother

\newbox\tmpbox
\newdimen\tmpdim
\def\narrow[#1]#2\par 
{%
\setbox\tmpbox\hbox{#2}%
\ifdim\wd\tmpbox>#1
\setbox\tmpbox\vbox{\hsize=#1 #2}%
\tmpdim=\ht\tmpbox%
\loop%
\setbox\tmpbox\vbox{\hsize=\wd\tmpbox \advance\hsize by -1pt
#2}%
\ifdim\ht\tmpbox=\tmpdim%
\relax%
\repeat%
\setbox\tmpbox\vbox{\hsize=\wd\tmpbox \advance\hsize by 1pt #2}%
\fi%
\box\tmpbox
}
\def\disp#1{
                $$
                \setbox\tmpbox\vbox{\narrow[\hsize]\noindent#1\par}
                \box\tmpbox
                $$
          }

\title{Strongly verbally closed groups
}
\author{Andrey M. Mazhuga
}
\date{\small
      Faculty of Mechanics and Mathematics\\
      Moscow State University\\
      Moscow 119991, Leninskie gory, MSU\\
      mazhuga.andrew@yandex.ru
}

\makeatletter
\def\@xfootnote[#1]{%
  \protected@xdef\@thefnmark{#1}%
  \@footnotemark\@footnotetext}
\makeatother

\begin{document}

\maketitle
\begin{otherlanguage}{english}

\begin{abstract}
It was recently proven that all free and many virtually free verbally closed subgroups are algebraically closed in \emph{any} group. We establish sufficient conditions for a group that is an extension of a free non-abelian group by a group satisfying a non-trivial law  to be algebraically closed in any group in which it is verbally closed. We apply these conditions to prove that the fundamental groups of all closed surfaces, except the Klein bottle, and almost all free products of groups satisfying a non-trivial law are algebraically closed in any group in which they are verbally closed.
\end{abstract}

\end{otherlanguage}

{\def\thefootnote{}\addtocounter{footnote}{-1}%
\footnote{%
The work of the author was supported by the Russian Foundation for Basic Research,
project no. 15-01-05823.} }

\section{Introduction}


A subgroup $H$ of a group $G$ is called \emph{verbally closed} (in $G$) \cite{[MR14]} (see also \cite{[Rom12]}, \cite{[RKh13]}, \cite{[KM17]}, \cite{[Mazh17]}) if any equation of the form $w(x_1,\dots,x_n)=h$, where $w(x_1,\dots,x_n)\in F_n(x_1,\dots,x_n)$ and $h\in H$, having a solution in $G$ has a solution in $H$ too.

A subgroup $H$ of a group $G$ is called \emph{algebraically closed} (in $G$) if any system of equations of the form $\left \{ w_1\left (x_1,\dots,x_n, H\right)=1,\dots,w_m\left ( x_1,\dots,x_n, H\right)=1 \right \}$, where $w_i\left (x_1,\dots,x_n,H\right) \in F_n(x_1,\dots,x_n)\ast H$, having a solution in $G$ has a solution in $H$ too.

A subgroup $H$ of a group $G$ is called \emph{retract} (of $G$) if $G$ is a semidirect product of a normal subgroup $N$ and  $H$ (i.e. $G=N\rtimes H$).

It is easy to see that any retract is an algebraically closed subgroup and any algebraically closed subgroup is a verbally closed subgroup. Thus, the question naturally arises: under what conditions on a subgroup $H$ and a group $G$ the revers implications hold. It is known (see, e.g., \cite{[KM17]}) that both of these reverse implications do not hold in general. The following was established in \cite{[MR14]}
{\it
\begin{enumerate}
    \item[\sffamily $\mathrm{R1)}$] if $G$ is finitely presented and $H$ is a finitely generated algebraically closed subgroup in $G$, then $H$ is a retract of $G$;

    \item[\sffamily $\mathrm{R2)}$] if $G$ is finitely generated over $H$\footnote[$\dagger$]{A group $G$ is \emph{finitely generated over $H$} if $G=\left \langle H,X \right \rangle$ for some finite subset $X\subseteq G$.} and $H$ is a equationally Noetherian\footnote[$\ddagger$]{A group $H$ is \emph{equationally Noetherian} if any system of equations with coefficients from $H$ and finitely many unknowns is equal to its finite subsystem.} algebraically closed subgroup in $G$, then $H$ is a retract of $G$.
\end{enumerate}
}
For the class of verbally closed subgroups, no similar structural descriptions in known. However, in finitely generated free groups (see \cite{[MR14]}) and finitely generated free nilpotent groups (see \cite{[RKh13]}) the situation is rather simple: verbally closed subgroups, algebraically closed subgroups and retracts are the same things. In article \cite{[KM17]}, the following theorem was established
\begin{ThKM}
Let $G$ be any group and let $H$ be its verbally closed virtually free infinite non-dihedral\footnote[$\dagger$]{for an infinite group, \emph{non-dihedral} means non-isomorphic to the free product of two groups of order two.} subgroup containing no infinite abelian noncyclic subgroups. Then
\begin{enumerate}
\item[\rm1)]
$H$ is algebraically closed in $G$;
\item[\rm2)]
if $G$ is finitely generated over $H$, then $H$ is a retract of $G$.
\end{enumerate}
\end{ThKM}
We call $H$ a \emph{strongly verbally closed group} if $H$ is an algebraically closed subgroup in any group containing $H$ as a verbally closed subgroup. Notice, that the first assertion of the given above theorem describes a certain class of strongly verbally closed groups (in particular, all non-trivial free groups belong to this class). Abelian groups is another class of strongly verbally closed groups (see Corollary \ref{CAg}). In this article we establish sufficient conditions for a group $H$ that is an extension of a free non-abelian group by a group satisfying a non-trivial law\footnote[$\ddagger$]{ Let $I(x_1,\dots,x_r)$ be an element of the free group $F_r(x_1,\dots,x_r)$ with a basis $x_1,\dots,x_r$, then we say that a \emph{group $G$ satisfies the law} $I$ if $I(g_1,\dots,g_r)=1$ in $G$ for all $g_1,\dots,g_r\in G$.} to be strongly verbally closed. These conditions can be used to establish strong verbal closedness of rather wide class of groups. For instance, is Section 2 we apply them to prove the following
{\it
\begin{enumerate}
 \item[\sffamily $\mathrm{C1)}$]
the fundamental groups of all closed surfaces, except the Klein bottle, are strongly verbally closed;
 \item[\sffamily $\mathrm{C2)}$]
if $H_i$, $i=1,\ \dots,\ n$, $n\ge 2$ are non-trivial groups satisfying a non-trivial law and $H=H_1\ast\cdots\ast H_n\not\simeq D_\infty$, then $H$ is strongly verbally closed.
\end{enumerate}
}
Notice, that any fundamental group $H$ of a closed surface is linear because it admits faithful representation in $\textrm{PSL}_2(\mathbb{R})$ (see, e.g., \cite{[CamKo07]}), and it is well-known that $\textrm{PSL}_2(\mathbb{R})$ is linear. Therefore, $H$ is equationally Noetherian \cite{[BMRe99]}. Thus, from C1) and R2) we have

\begin{myth} If $H$ is the fundamental group of a closed surface, which is not the Klein bottle, then $H$ is a retract of any finitely generated over $H$ group containing $H$ as a verbally closed subgroup.
\end{myth}

Before proceeding to the formulation of the main result, we introduce some notation.
The symbol $\left \langle x \right \rangle_\infty$ means an infinite cyclic group generated by $x$. The symbols $I_A(H)$ and $I_A(\left \langle x \right \rangle_\infty)$ mean the verbal subgroups of $H$ and $\left \langle x \right \rangle_\infty$ generated by a word $I_A$\footnote[$\dagger\dagger$]{ Let $I(x_1,\dots,x_r)$ be an element of the free group $F_r(x_1,\dots,x_r)$ with a basis $x_1,\dots,x_r$, then the \emph{verbal subgroup} of a group $G$ determined by $I$, $I(G)$, is $I(G)=\left \langle I(g_1,\dots,g_r)|\ g_1,\dots,g_r\in G \right \rangle.$}, respectively. The symbol $I_A(H)^{\langle x  \rangle_\infty\ast I_A(H)}$ means the normal closure of $I_A(H)$ in $\langle x  \rangle_\infty\ast I_A(H)$. The symbols $Z(H)$ and $H'$ mean the center and the derived subgroup of $H$, respectively.


\begin{Main Theorem}\label{MTh} If a group $H$ such that there exists
\begin{enumerate}
 \item[\sffamily $\mathrm{T1)}$] a short exact sequence of groups $1\to F\to H\to A\to 1$ such that $F$ is a free non-abelian group and $A$ is a group satisfying a non-trivial law $I_A$;

 \item[\sffamily $\mathrm{T2)}$] a divisible subgroup $Q$ of $H$ such that $Q\subseteq Z(H)$ and $Q\cap H'=\{1\}$;

 \item[\sffamily $\mathrm{T3)}$] a generation set $U=U^{-1}$ of $H$ such that for any $u\in U$ there are elements $E_{u,k}\left (x,I_A\left ( H \right )  \right )\in  I_A(\langle x  \rangle_\infty)\cdot I_A(H)^{\langle x  \rangle_\infty\ast I_A(H)}$, $k=1,\ \dots,\ n_u$ such that the system of equations $$\left \{ E_{u,k}\left (x,I_A\left ( H \right )  \right )=E_{u,k}\left (u,I_A\left ( H \right )  \right )|\ k=1,\ \dots,\ n_u \right \}$$ has a unique solution in $H/Q$,
\end{enumerate} then $H$ is strongly verbally closed.
\end{Main Theorem}
Hence, to use this theorem one has to choose at one's discretion: 1) a generating set $U=U^{-1}$ of $H$, 2) a non-trivial law, $I_A$, in $A$, 3) a central divisible subgroup $Q$ of $H$ such that $Q\cap H'=\{1\}$ (it can be shown that in this case $Q$ is a direct factor of $H$, however, we will neither need or use this fact), 4) a set of words $\left \{ E_{u,k}\left (x,I_A\left ( H \right )  \right )|\ k=1,\ \dots,\ n_u \right \}\subset I_A(\langle x  \rangle_\infty)\cdot I_A(H)^{\langle x  \rangle_\infty\ast I_A(H)}$ for each $u\in U$ (note that these sets may differ for the different elements $u\in U$).

Now we want to draw the reader's attention to the fact that at the definition of a strongly verbally closed group $H$ verbal closedness of $H$ (in some group $G$) is given as the assumption. We do not consider the question in which groups a given group $H$ can be verbally closed, and we restrict ourselves to the following remark. It follows from Example \ref{Ex2} that any free group is strongly verbally closed. On the other hand

{\Observation If a group $G$ has a generating set consisting of finite-order elements, then $G$ does not have a verbally closed non-trivial free subgroup.}
\begin{Proof} Let $F_r$ be a verbally closed free subgroup of $G$ of rank $r\ge 1$ and let $S$ be a generating set of $G$ consisting of finite-order elements. If $r=1$, then $F_r=\left \langle f \right \rangle_\infty$ for some $f\in G\setminus \{1\}$. Let $s_{1}s_{2}\cdots s_{m}$ be a decomposition of $f$ into the elements $s_i\in S.$
Then, there is a natural number $p\ge 2$ which is relatively prime to the orders, $\textrm{ord}_G(s_{i})$, of all the elements $s_{i}$ of this decomposition. Consider the equation: $$x_{1}^px_{2}^p\cdots x_{m}^p=f.$$
This equation has a solution in $G$ (it is easy to see that $x_i=s_{i}^{m_i}$, where $pm_i\equiv 1\pmod{\textrm{ord}_G(s_{i})}$ is a solution of the equation), but it has no solution in $\left \langle f \right \rangle_\infty$ (since the element $f$ is not a proper power in $\left \langle f \right \rangle_\infty$). Thus, this contradicts the verbal closedness of $H$ in $G$.

If $r\ge 2$, then $F_r = \left \langle f \right \rangle_\infty\ast F_{r-1}$, for some $f\in G\setminus \{1\}$ and a free subgroup $F_{r-1}$ of $G$. The subgroup $\left \langle f \right \rangle_\infty$ is verbally closed in $F_r$ (since $\left \langle f \right \rangle_\infty$ is a free factor of $F_r$, and a free factor is a retract), therefore (by the transitivity of verbal closedness), it is verbally closed in $G$ and we again get a contradiction.
\end{Proof}

In Section 2 we discuss examples and corollaries showing how the Main Theorem can be applied. Section 3 contains the proof of the Main Theorem. Our argument is a refined and more sophisticated version of the proof of the Main Theorem of \cite{[KM17]}, but it also based on the use of Lee words \cite{[Lee02]}.

\hfill \break

Let us fix the notations. If $k\in \mathbb{Z}$, $x$ and $y$ are elements of a group, then $x^y$, $x^{ky}$, $x^{-y}$ and $x^{-ky}$ denote $y^{-1}xy$, $y^{-1}x^ky$, $y^{-1}x^{-1}y$ and $y^{-1}x^{-k}y$, respectively. The derived subgroup and the center of a group $G$ are denoted by $G'$ and $Z(G)$, respectively. The commutator, $[x,y]$, of two elements, $x$ and $y$, of a group we define as $x^{-1}y^{-1}xy.$ If $X$ is a subset of a group $G$, then $\left | X \right |$, $\left \langle X \right \rangle$ and $C_G(X)$ mean the cardinality of $X$, the subgroup generated by $X$ and the centralizer of $X$ in $G$, respectively. Notations $\left \langle x \right \rangle_n$, $n\in\mathbb{N}$ and $\left \langle x \right \rangle_\infty$ mean a cyclic group of order $n$ and an infinite cyclic group generated by an element $x$, respectively. The index of a subgroup $H$ of a group $G$ is denoted by $(G:H)$. If $A$ and $B$ are subgroups of a group $G$, then $A\cdot B=\left \langle ab|a\in A,\ b\in B \right \rangle$. The symbol $A^B$ means the normal closure of a subgroup $A$ in a group $B$. The symbols $A\ast B$, $A\underset{C}{\ast}B$, $A\rtimes B$ и $A\times B$ mean the free product of groups $A$ and $B$, a free product of groups $A$ and $B$ with the amalgamated subgroup $C$, a semidirect product of groups $A$ and $B$ and the direct product of groups $A$ and $B$, respectively. $F_n(x_1,\dots,x_n)$ or $F_n$ is the free group of rank $n$ (with a basis $x_1,\dots,x_n$). Sometimes instead of a notation $x_1,\dots,x_n$ we use an abridged notation $\underline{x}$ rewriting, for example, $w(x_1,\dots,x_n)$ as $w(\underline{x})$ or $F_n(x_1,\dots,x_n)$ as $F(\underline{x})$.

\section{Examples and Corollaries}

The following example has no particular theoretical interest and was included to demonstrate how the Main Theorem can be applied to a concrete group.

{\Example\label{A0505FE} The group $H = F_2(b,c)\rtimes\left \langle a \right \rangle_2$, where $F_2(b,c)$ is the free group freely generated by $b,c$ and $b^a = b^{-1}$, $c^a=c$ is strongly verbally closed.
}

\begin{Proof} Consider the following short exact sequence of groups
$$
1\to F_2(b,c)\to F_2(b,c)\rtimes\left \langle a \right \rangle_2\to \left \langle a \right \rangle_2\to 1
$$
and choose: 1) the generating set $U=U^{-1}=\left \{ (b,a),\ (c,a),\ (c^{-1},a),\ (1,a) \right \}$, 2) the law $I_A = t^2$ and 3) the group $Q=\{1\}$. Next, we choose for each $u\in U$ the words $E_{u,k}$ and verify that the conditions of the Main Theorem hold.

For $(b,a)$ we consider the following words: $E_{(b,a),1}=x^2,\ E_{(b,a),2}=\left [x,(b^2,1)  \right ],\ E_{(b,a),3}=\left [x,(c^2,1)  \right ]\in I_A(\langle x  \rangle_\infty)\cdot I_A(H)^{\langle x  \rangle_\infty\ast I_A(H)}$ and show that the system
$$
\left\{ x^2=(b,a)^2=(1,1),\ \left [x,(b^2,1)  \right ]=\left [(b,a),(b^2,1)  \right ]=(b^{4},1),\ \left [x,(c^2,1)  \right ]=\left [(b,a),(c^2,1)  \right ]=(\left [b,c^2  \right ],1)\right\}
$$
has the unique solution, $x=(b,a)$, in $H/Q=H$. It is straightforward to verify that the solutions of the equation $x^2=(1,1)$ in $H$ are either of the form $x=(w_1b^kw_2,a)$, where $k\in\mathbb{Z}$ and $w_1^a=w_2^{-1}$ or $x=(1,1)$. It is clear that $x=(1,1)$ is not a solution of the second equation. Substituting $x=(w_1b^kw_2,a)$ in $\left [x,(b^2,1)  \right ]=(b^{4},1)$ we have $w_1b^kw_2b^{2}w_2^{-1}b^{-k}w_1^{-1}b^{2}=b^{4}$ in $F_2(b,c)$, that is $b^{2w_2^{-1}b^{-k}w_1^{-1}}=b^{2}$ in $F(b,c)$, whence $w_1b^{k}w_2=b^m$ for some $m\in\mathbb{Z}$. Since we can assume that the word $w_1b^{k}w_2$ is reduced and that the last letter of $w_1$ and the first letter of $w_2$ is not $b$, the last equality means that $w_1=w_2=1$ (and $k=m$). Thus, the solutions of the first two equations of the system are of the form  $x=(b^k,a)$, $k\in\mathbb{Z}$. Substituting $x=(b^k,a)$ in the third equation we obtain the equality $\left [b^k,c^2  \right ]=\left [b,c^2  \right ]$ in $F_2(b,c)$, from which it follows that $k=1$.

A similar verification shows that one can choose the words $E_{(c,a),1}=E_{(c^{-1},a),1}=x^2,\ E_{(c,a),2}=E_{(c^{-1},a),2}=\left [x,(b^2,1)  \right ]$ for $(c,a)$ and $(c^{-1},a)$, and the words $E_{(1,a),1}=x^2,$ $E_{(1,a),2}=\left [x,(b^2,1)  \right ],$ $E_{(1,a),3}=\left [x,(c^2,1)  \right ]$ for $(1,a)$.
\end{Proof}

The following two examples are the special cases of assertion 1) of the stated above theorem of \cite{[KM17]}.

{\Example\label{Ex2} A free group is strongly verbally closed.
}
\begin{Proof} Let $H$ be a free group. If $H$ is abelian, then $H$ is strongly verbally closed in accordance with Corollary \ref{CAg}. Let $H$ be a non-abelian free group. Consider the following short exact sequence of groups
$$1\to H\to H\to 1\to 1$$
 and choose $U= H,$ $I_A=t$, $Q=\{1\}$ and $E_{u,1}=x \in I_A(\langle x  \rangle_\infty)\cdot I_A(H)^{\langle x  \rangle_\infty\ast I_A(H)}$ for each $u\in U$. Now, the verification of the conditions of the Main Theorem is trivial.
\end{Proof}

{\Example Virtually free non-virtually cyclic group with the unique extraction of roots of infinite-order elements\footnote[$\dagger$]{that is, if $h_1,h_2$ are the infinite-order elements and $h_1^k=h_2^k$, $k\in\mathbb{N}$, then $h_1=h_2$.} is strongly verbally closed.
}

\begin{Proof} If $H$ is a virtually free non-virtually cyclic group, then it is clear that there exists a short exact sequens of groups of the form
$$
1\to F\to H\to A\to 1,
$$
where $F$ is a non-abelian free group, and $A$ is a finite group. Let $U$ be the set of all infinite-order elements of $H$, $I_A=t^{\left | A \right |}$, $Q=\{1\}$, and $E_{u,1} =  x^{\left | A \right |}\in I_A(\langle x  \rangle_\infty)\cdot I_A(H)^{\langle x  \rangle_\infty\ast I_A(H)}$ for each $u\in U$.

In \cite{[KM17]}, the following was established
\disp{\hfuzz15pt\sl
in a virtually free group which is not virtually cyclic, any element decomposes into a product of two infinite-order elements.
}
It follows immediately from this statement that $U$ is a generating set of $H$. The uniqueness of a solution of an equation $x^{\left | A \right |}=u^{\left | A \right |}$, $u\in U$ in $H/Q=H$ is assumed in the hypothesis of the example. The verification of the remaining conditions of the Main Theorem is trivial.
\end{Proof}

{\Corollary\label{Cor3} Let $1\rightarrow F\rightarrow H \rightarrow  A\rightarrow 1$ be a short exact sequence of groups such that $F$ is a non-abelian free group, and $A$ is a group satisfying a non-trivial law  $I_A$. If $Z(H)$ is a divisible group, $Z(H)\cap H'=\{1\}$ and there exist elements $f_1,\ \dots,\ f_m\in I_A(H)$ such that $C_H\left ( f_1,\dots,f_m \right )=Z\left ( H \right )$, then $H$ is strongly verbally closed.
}

\begin{Proof} Let us put $U=H$, $Q = Z(H)$ and $E_{u,k}=\left [ f_k,x \right ]\in I_A(\langle x  \rangle_\infty)\cdot I_A(H)^{\langle x  \rangle_\infty\ast I_A(H)}$, $k=1,\ \dots,\ m$ for each $u\in U$.

To verify the conditions of the Main Theorem we only need to check that for $u\in U$ the system
$$
\left \{ \left [ f_k,x \right ]=\left [ f_k,u \right ] |\ k=1,\ \dots,\ m \right \}
$$
has a unique solution in $H/Q$. It is clear that the system has a solution, $x=uQ$, in $H/Q$. Let $\widehat{x}Q$ be a solution of this system in $H/Q$, then (since $Q=Z(H)$, and $Z(H)\cap H'=\{1\}$) for all $k$ we have $f_k^{-1}\widehat{x}^{-1}f_k\widehat{x}=f_k^{-1}u^{-1}f_ku$ in $H$, thus $\widehat{x}u^{-1}\in C_H(f_k)$, $k=1,\ \dots,\ m$. That is, $\widehat{x}u^{-1}\in\bigcap_{k=1}^{m}C_H(f_k)=C_H(f_1,\dots,f_m)=Q$, therefore $\widehat{x}Q=uQ$ in $H/Q$.

\end{Proof}

Further we will need the following technical lemma.

\begin{Lemma}\label{LemWSF} Let $1\rightarrow F\rightarrow H \rightarrow  A\rightarrow 1$ be a short exact sequence of groups such that $F$ is a non-abelian free group, and $A$ is a group satisfying a non-trivial law  $I_{A}$. Then $I_{A}(H)$ is a non-abelian free subgroup of $H$.
\end{Lemma}
\begin{Proof} The group $I_{A}(H)$ is free as a subgroup of a free group $F$. Subgroup $I_{A}(F)$ in normal (even fully invariant) in $F$ as a verbal subgroup and $I_{A}(F)\ne\{1\}$ since $I_{A}$ is a non-trivial word. Therefore, if $I_{A}(H)$ is abelian, then $I_{A}(F)\simeq\mathbb{Z}$. It can be easily deduce from the Nielsen-Schreier Theorem that a non-abelian free group cannot have a normal infinite cyclic subgroup. Thus, $I_{A}(H)$ is a non-abelian free group.
\end{Proof}

{\Corollary\label{Cor4} Let $H = H_1\ast\cdots\ast H_n\not\simeq D_\infty$, $n\ge2$ be a free product of non-trivial groups $H_i$ satisfying a non-trivial law, then $H$ is strongly verbally closed.
}
\begin{Proof} Consider the following short exact sequence of groups:
$$
1\to F\to H\to A\to 1,
$$
where $A=H_1\times\cdots\times H_n$ and $F$ is the kernel of the natural epimorphism  $H_1\ast\cdots\ast H_n\twoheadrightarrow A$ (that is, $F$ is the cartesian subgroup). It is well-known that the cartesian subgroup is a free group (and it is clear that if $H\not\simeq D_\infty$ and $n\ge 2$, then $F$ is non-abelian free).

Let $[x_1,x_2,\dots,x_{n-1},x_n]=\left[ \left[ \left [\dots\left [ \left [ x_1,x_2 \right ],x_3 \right ],\dots\right],x_{n-1}\right],x_{n}\right]$ be the left-normed commutator of $x_1,x_2,\dots,x_{n-1},x_n$, and let $W_i$ be a non-trivial law in $H_i$ (we assume that for the different indices the words $W_i$ do not contain the same letters). It is easy to see that in this case $A$ satisfies the non-trivial law $I_A = \left[W_1,W_2,\dots,W_{n-1},W_n\right]$. We shall verify the conditions of Corollary \ref{Cor3}.

The group $B=I_A\left ( H \right )$ is non-abelian free by Lemma \ref{LemWSF}. Therefore, there exist elements $f_1,f_2\in B$ such that $[f_1,f_2]\ne1$. It is well-known that the cartesian subgroup trivially intersects the free factors $H_i$, therefore (since $B$ is a subgroup of the cartesian subgroup $F$), the elements $f_1,f_2$ cannot belong to a subgroup of the form $H_i^g$, $g\in H$. If $C_H(f_1)\cap C_H(f_2)\ne\{1\}$, then the center of the group $\left \langle f_1,f_2,z \right \rangle$, where $z\in C_H(f_1)\cap C_H(f_2)\setminus \{1\}$ is non-trivial, therefore the group $\left \langle f_1,f_2,z \right \rangle$ does not decompose as a proper free product. But by the Kuro\u{s} Subgroup Theorem (see. e.g., \cite{[Rob95]}) this would imply that $\left \langle f_1,f_2,z \right \rangle$ is a subgroup of $H_i^g$, $g\in H$. Thus, $C_H(f_1)\cap C_H(f_2)=C_H(f_1,f_2)=1.$


\end{Proof}

 Sice the fundamental group of a connected non-closed surface is free (see, e.g., \cite{[Stillwell]})), by Example \ref{Ex2}, it is strongly verbally closed. It is well-known (see, e.g., \cite{[LS77]}) that the fundamental group of a closed surface does not decompose as a proper free product, nevertheless

{\Corollary\label{FSC} The fundamental group of a closed surface, except the Klein bottle, is strongly verbally closed.
}

Further in the text of the article a surface means a closed (i.e. connected, compact, without boundary) surface.

It is well-known that the fundamental group of an orientable surface of genus $g\ge0$ (we denote such a surface as $S_g^+$) is of the type:
$$
\pi_1(S_g^+)=\left \langle x_1,y_1,\dots,x_g,y_g|\ [x_1,y_1]\cdots[x_g,y_g] \right \rangle,
$$
and that the fundamental group of a non-orientable surface of genus $g\ge1$ (we denote such a surface as $S_g^-$) is of the type:
$$
\pi_1(S_g^-)=\left \langle x_1,\dots,x_g|\ x_1^2\cdots x_g^2 \right \rangle.
$$
The following technical lemma is rather well-known.
\begin{Lemma}\label{FSLemma} Let $\pi_1(S)$ be the fundamental group of a surface $S$ whose Euler characteristic\footnote[$\dagger$]{The Euler characteristic, $\chi(S_g^+)$, of $S_g^+$ is $\chi(S_g^+)=2-2g$ and the Euler characteristic, $\chi(S_g^-)$, of $S_g^-$ is $\chi(S_g^-)=2-g$.} is $\chi(S)$.

    \begin{enumerate}
\item[\sffamily $\mathrm{S1)}$] If $\chi(S)< 0$ and elements $g_1,g_2\in \pi_1(S)\setminus\! \left \{ 1 \right \}$ commute, then $\left \langle g_1,g_2 \right \rangle\simeq\mathbb{Z}$.
\item[\sffamily $\mathrm{S2)}$] If $\chi(S)< 0$, then $\pi_1(S)'$ is a non-abelian free group.
    \end{enumerate}
\end{Lemma}
\begin{Proof} For a proof of S1) see, for example, \cite{[Jaco70]}.

It is well-known (see, e.g., \cite{[Jaco70]}) that $\pi_1(S)'$ is a free group. If $\chi(S)< 0$, then $g\ge 2$ for $S_g^+$ and $g\ge 3$ for $S_g^-$. Therefore, the following epimorphisms are naturally defined
\begin{align*}
\pi_1(S_g^+)=\left \langle x_1,y_1,\dots,x_g,y_g|\ [x_1,y_1]\cdots[x_g,y_g] \right \rangle &\twoheadrightarrow \left \langle x_1,y_1,x_2,y_2\ |\ [x_1,y_1][x_2,y_2] \right \rangle\simeq F_2\underset{\mathbb{Z}}{ \ast }F_2;\\
\pi_1(S_g^-)=\left \langle x_1,\dots,x_g|\ x_1^2\cdots x_g^2 \right \rangle &\twoheadrightarrow \left \langle x_1,x_2,x_3\ |\ x_1^2,x_2^2,x_3^2 \right \rangle\simeq \mathbb{Z}_2\ast\mathbb{Z}_2\ast\mathbb{Z}_2.
\end{align*}
Since the commutator subgroups of the images of these epimorphisms are non-abelian groups, it follows that the commutator subgroups, $\pi_1(S_g^+)'$, $g\ge 2$ and $\pi_1(S_g^-)'$, $g\ge 3$, are also non-abelian.
\end{Proof}

\begin{Proof} It is clear that the groups $\pi_1(S_0^+)$, $\pi_1(S_1^+)$ and $\pi_1(S_1^-)$ are abelian, therefore, by Corollary \ref{CAg}, these groups are strongly verbally closed. The group $\pi_1(S_2^-)$ is excluded by the hypothesis of the corollary. Thus, it remains to consider the surfaces whose Euler characteristic is negative.

For $H=\pi_1(S)$, $\chi(S)<0$ consider the following short exact sequence of groups
$$
1\to H'\to H \to A\to 1,
$$
where $A=H/H'$. By assertion S2) of Lemma \ref{FSLemma}, $H'$ is a non-abelian free group. Assuming that $I_A = [t_1,t_2]$, we verify the conditions of Corollary \ref{Cor3}.

By Lemma \ref{LemWSF}, $B=I_A(H)$ is a non-abelian free subgroup of $H$ and, therefore, there are elements $f_1,f_2\in B$ such that $[f_1,f_2]\ne1$. Suppose that there exists an element, say $z$, such that $z\in C_H(f_1)\cap C_H(f_2)\setminus \{1\}.$ Then, by assertion S1) of Lemma \ref{FSLemma}, there are elements $u_1, u_2\in H$ such that $f_1=u_1^{k_1}$, $z=u_1^{m_1}$ and $f_2=u_2^{k_2}$, $z=u_2^{m_2}$ for some $k_1,k_2,m_1,m_2\in\mathbb{Z}\setminus\! \{0\}$. Therefore, we have $f_1^{k_2m_1}=u_1^{k_1k_2m_1}=z^{k_1k_2}=u_2^{k_1k_2m_2}=f_2^{k_1m_2}$. Since $f_1$ and $f_2$ are elements of a free group, the equality $f_1^{k_2m_1}=f_2^{k_1m_2}$ means that $f_1$ and $f_2$ commute, but this contradicts our choice. Whence, $C_H(f_1,f_2)=C_H(f_1)\cap C_H(f_2)=\{1\}$.
\end{Proof}

\section{Proof of the Main Theorem}

The following lemma is well known.
\begin{Lemma}\label{A0505WA}
If a subgroup $H$ of a group $G$ is such that any finite system
of equations of the
form
\begin{equation}\label{A0505SES}
\left \{ w_i(x_1,\dots,x_n)=h_i\ |\ i=1,\ \dots,\ m \right \},
\end{equation}
where $w_i(x_1,\dots,x_n)\in F_n(x_1,\dots,x_n)$, $h_i\in H$,
having a solution in $G$
has a solution in $H$ too, then $H$ is algebraically closed.
\end{Lemma}

\begin{Proof}
Just denote the coefficients by new letters and interpret them as
variables. For example, the solvability of the equation
$xyh_1[x^{\the\year},h_2]y^{-1}=1$
is equivalent to the solvability of the system
$$
\{xyz[x^{\the\year},t]y^{-1}=1,\
z=h_1,\
t=h_2\}.
$$
\end{Proof}

Now, recall that any integer matrix can be reduced to a diagonal matrix by integer elementary transformations. This means that any finite system of equations of the form (\ref{A0505SES}) can be reduced to a system of the form
\begin{equation}\label{A0505SEM}
\left \{ x_i^{m_i}u_i(x_1,\dots,x_n)=h_i,\ u_j(x_1,\dots,x_n)=h_j \ |\ i=1,\ \dots,\ l,\ j=l+1,\ \dots,\ m  \right \},
\end{equation}
where $u_i(x_1,\dots,x_n),u_j(x_1,\dots,x_n)\in F(x_1,\dots,x_n)'$, $m_i>0$, $0\le l\le n$, $h_i,h_j\in H$, by means of a finite sequence of transformations of the form $w_i\to w_iw_j^{\pm 1}$ and $x_i\to x_ix_j^{\pm 1}$. It follows from the form of these transformations that system (\ref{A0505SES}) has a solution in a group $G$ (containing $H$ as a subgroup) if and only if system (\ref{A0505SEM}) has a solution in $G$. Therefore, the following lemma holds.

\begin{Lemma}\label{LemWA2}
If a subgroup $H$ of a group $G$ is such that any finite system of the form (\ref{A0505SEM}) having a solution in $G$ has a solution in $H$ too, then $H$ in algebraically closed.
\end{Lemma}

{\Corollary\label{CAg} An abelian group is strongly verbally closed.}

\begin{Proof} Let $H$ be a verbally closed abelian subgroup in a group $G$. Suppose that a system of the form (\ref{A0505SEM}) has a solution in $G$. Since $H$ is a verbally closed subgroup in $G$, each equation of this system has a solution in $H$. For each $i=1,\ \dots,\ l$ let a tuple $\left ( a_{i,1},a_{i,2},\dots,a_{i,n} \right )$ be a solution of the equation $x_i^{m_i}u_i(x_1,\dots,x_n)=h_i$ in $H$. Since $u_i(x_1,\dots,x_n),u_j(x_1,\dots,x_n)\in F(x_1,\dots,x_n)'$ and $H$ is abelian, we have $u_i(g_1,\dots,g_n)=1$, $u_j(g_1,\dots,g_n)=1$ for any $g_1,\dots,g_n\in H$, thus $h_j=1$, $j=l+1,\ \dots,\ m$. Now it is easy to see that $x_i=a_{i,i}$, $i=1,\ \dots,\ l$ and $x_i=1$, $i=l+1,\ \dots, m$ is a solution of (\ref{A0505SEM}) in $H$.
\end{Proof}

Recall that a \emph{Lee word}
in $m$ variables
for the free
group of rank $r$
is an element
$L_m(z_1,\dots,z_m)$ of the free group of rank $m$ such that
\goodbreak
\begin{itemize}
\item[L1)]
if $L_m(v_1,\dots,v_m)=L_m(v_1',\dots,v_m')\ne1$ in $F_r$, then
$v_i'\in F_r$ are obtained from $v_i\in F_r$ by
simultaneous conjugation, i.e.,
there exists $s\in F_r$ such that
$v_i'=v_i^s$ for all $i=1,\ \dots,\ m$;

\nobreak
\item[L2)]
$L_m(v_1,\dots,v_m)=1$ if and only if the elements
$v_1,\dots,v_m$ of $F_r$ generate a cyclic subgroup.

\end{itemize}

\goodbreak

In \cite{[Lee02]}, such words were constructed for all
integers $r,m\ge2$.
Actually, it is easy to see that Lee's result implies
the existence of
a \emph{universal Lee word} in $m$ variables.

\begin{UTW}
For any positive integer $m$, there exists
an element
$L_m(z_1,\dots,z_m)\in F_m$
such that
properties {\rm L1)} and {\rm L2)} hold in
all free groups
$F_r$ and even in $F_\infty$.
\end{UTW}

\begin{Proof}
This assertion follows immediately from Lee's result and the following
simple fact:
\disp{\sl
$F_\infty$ embeds into $F_2$ as a \emph{malnormal} subgroup,
}%
i.e. a subgroup
$S\subset F_2$
such that $S^f\cap S=\{1\}$ for all $f\in F_2\setminus S$.
This fact
follows, e.g., from a result of~\cite{[Wise01]}:
\disp{\sl
\hfuzz 6pt
in a free group, any set
satisfying small-cancellation condition $C(5)$ freely
generates a malnormal subgroup.
}
Thus, a Lee word for $F_2$ is universal, i.e. it is suitable also
for $F_\infty$.
\end{Proof}

The proof of the following lemma is based on the use of universal Lee words.

{\Lemma\label{SmL} Let $I(H)$ and $I\left ( F(\underline{x}) \right )$ be the verbal subgroups of $H$ and $F(\underline{x})$ generated by a word $I$, respectively. If $I(H)$ is a non-abelian free group and $H$ is a verbally closed subgroup in a group $G$, then any finite system of equations of the form
\begin{equation}\label{Nth01}
\left \{ w_i(\underline{x},I(H))=1|\ i=1,\ \dots,\ m \right \},
\end{equation}
where $w_i(\underline{x},I(H))\in I\left( F(\underline{x}) \right )\cdot I(H)^{F(\underline{x})\ast I(H)}$, having a solution in $G$ has a solution in $H$ too.
}
\begin{Proof} Let $S$ be a system of the form (\ref{Nth01}), having a solution in $G$.

By the hypothesis of the lemma, the left-hand sides, $w_i(\underline{x},I(H))$, of the equations of $S$ have decompositions of the form
\begin{equation}\label{Neq01}
w_i(\underline{x},I(H))=a_i\prod_{j=1}^{m_i}b_{i,j}^{c_{i,j}},
\end{equation}
 where $a_i\in I\left( F(\underline{x}) \right )$, $b_{i,j}\in I(H)$ and $c_{i,j}\in F(\underline{x})\ast I(H)$. We fix such a decomposition for the left-hand side of each equation of $S$. Let $M_{i,j}$ be the set of all the elements of $I(H)$ which are present in the normal form (for $F(\underline{x})\ast I(H)$) of the element $c_{i,j}$. We define the set $M$ as the union of all the sets $M_{i,j}$ and all the elements $b_{i,j}$. Since $I(H)$ is a non-abelian free group, there exist elements $h',h''\in I(H)$ such that $\left \langle M,h',h'' \right \rangle\not\simeq\mathbb{Z}$. If $\left \langle M \right \rangle\simeq\mathbb{Z}$ or $M=\varnothing$, we add the elements $h'$, $h''$ to the set $M$.

Since $M\subset I(H)$, for each $h_\beta\in M$ there is a word, $t_\beta(\underline{y}_{\beta})\in I(F(\underline{y}_{\beta}))$, such that $t_\beta(\underline{g}_{\beta}) = h_\beta$ for some elements $\underline{g}_{\beta}$ of $H$ (we assume that for the different indices the tuples $\underline{y}_{\beta}$ do not contain the same letters).
For each element $h_\beta\in M$ we add to the system $S$ the equation $t_\beta(\underline{y}_{\beta})=h_\beta$, and each occurrence of the constant $h_\beta$ in the fixed above decompositions (\ref{Neq01}) of the left-hand sides of the equations of the system $S$ we replace with the word $t_\beta(\underline{y}_{\beta}).$ We denote this new system by $S^1$. Note that the left-hand sides of the equations of $S^1$ do not contain constants from $H$; their right-hand sides are the constants from $I(H)$; the system $S^1$ has a solution in $G$ if and only if the system $S$ has a solution in $G$.

Let us denote the equations of the system $S^1$ by $u_\alpha(\underline{x},\underline{y})=f_\alpha$, where $\underline{y}$ is a concatenation of all the tuples $\underline{y}_{\beta}$\footnote[$\dagger$]{if e.g., $(y_1,y_2)$ and $(y_3,y_4)$ are tuples, then we call the tuples of the form $(y_1,y_2,y_3,y_4)$ or $(y_3,y_4,y_1,y_2)$  their concatenation.}. That is
$$
S^1=\left \{ u_\alpha(\underline{x},\underline{y})=f_\alpha|\ \alpha=1,\ \dots,\ N \right \},
$$ where $N=m+|M|$.
Let $L_N(z_1,\dots,z_N$) be a universal Lee word in $N$ variables. Consider the equation
$$
L_N(u_1(\underline{x},\underline{y}),\dots,u_N(\underline{x},\underline{y}))=L_N(f_1,\dots,f_N).
$$
This equation has a solution in $G$ by construction (just take the following: a solution of $S$ as $\underline{x}$ and $\underline{g}_{\beta}$ as $\underline{y}_{\beta}$). Thus, since $H$ is a verbally closed subgroup in $G$ and $L_N(f_1,\dots,f_N)\in I(H)\subseteq H$, the equation has a solution, say $\underline{\hat{x}}$, $\underline{\hat{y}}$, in $H$.

By construction a word $u_\alpha(\underline{x},\underline{y})$ is either of the form
\begin{equation}\label{Neq03}
w_i(\underline{x},\underline{y})\in I\left ( F(\underline{x})\right )\cdot I(F(\underline{y}))^{F(\underline{x})\ast I(F(\underline{y}))}
\end{equation}
or of the form
\begin{equation}\label{Neq04}
t_\beta(\underline{y}_{\beta})\in I(F(\underline{y})).
\end{equation}
It is clear from (\ref{Neq03}) and (\ref{Neq04}) that for any values of $\underline{x}$, $\underline{y}$ in $H$ the value of $u_\alpha(\underline{x},\underline{y})$ is an element of the free group $I(H)$.


 Since $I(H)$ is a free group, $f_\alpha, u_\alpha(\underline{\hat{x}},\underline{\hat{y}})\in I(H)$ for all $\alpha=1,\ \dots,\ N$ and $\left \langle f_1,\dots,f_N \right \rangle\not\simeq\mathbb{Z}$ (indeed, by construction $M\subseteq \{f_1,\dots,f_N\}$ and $\left \langle M \right \rangle\not\simeq\mathbb{Z}$), Property L2) of Lee words implies that $L_N(f_1,\dots,f_N)\ne1.$ Then, according to Property L1) of Lee words, there exists $s\in I(H)$ such that
$$\hat{S}^1=\left \{ u_\alpha(\underline{\hat{x}},\underline{\hat{y}})=f_\alpha^s|\ \alpha=1,\ \dots,\ N \right \}.$$
That is, $\underline{\hat{x}}^{s^{-1}}$, $\underline{\hat{y}}^{s^{-1}}$ is a solution to the system $S^1$ in $H$. Recall how the system $S^1$ was obtained from the system $S$. Now we carry out the inverse transformation, but with the set of equalities $\hat{S}^1$. More precisely, using the equalities of the form $t_\beta(\underline{\hat{y}}_{\beta})=h_\beta^s$ of $\hat{S}^1$, in the left-hand sides of all the equalities of the form $w_i(\underline{\hat{x}},\underline{\hat{y}})=1^s$ of $\hat{S}^1$ we replace (according to the fixed above decompositions (\ref{Neq01}) of the left-hand sides of the equations of $S$) all the elements $t_\beta(\underline{\hat{y}}_{\beta})$ with the elements $h_\beta^s$. This leads us to the following set of equalities
$$
\hat{S}=\left \{ w_i(\underline{\hat{x}},I(H)^s)=1|\ i=1,\ \dots,\ m \right \}.
$$
Whence, $\underline{\hat{x}}^{s^{-1}}$ is a solution to the system $S$ in $H$ as required.
\end{Proof}

Now we are ready to prove the Main Theorem.

\begin{ProofMT}
Let $H$ be a verbally closed subgroup in a group $G$, then (by Lemma  \ref{LemWA2}) it suffices to show that a system $S$ of the form (\ref{A0505SEM}) having a solution in $G$ has a solution in $H$ too.

To each equation $w_i\left ( \underline{x} \right )=h_i$ of the system $S$ we associate the system

\begin{equation}\label{SEm1}
S_i=\left \{ w_i\left ( \underline{x} \right )y_{i,t_i}^{-1}\cdots y_{i,1}^{-1}=u_{i,0},\ y_{i,1}=u_{i,1},\ \dots,\ y_{i,t_i}=u_{i,t_i} \right \},
\end{equation}
where $h_i = u_{i,0}u_{i,1}\dots u_{i,t_i},$ $u_{i,j}\in U$ and $y_{i,j}$ are the new unknowns. Consider the system
$$
S^1=\left \{ S_i|\ i=1,\ \dots,\ m\right\},
$$
consisting of all the equations of all the systems $S_i.$ It is clear that the system $S^1$ has a solution in $G$ if and only if the system $S$ has a solution in $G$.

Using the elements $E_{u_{i,j},k}\left (x,I_A\left ( H \right )  \right )$, $k=1,\ \dots,\ n_{u_{i,j}}$ from Condition T3) of the Main Theorem, to each equation $w_{i,j}\left ( \underline{x},\underline{y} \right )=u_{i,j}$ (where $w_{i,0}\left ( \underline{x},\underline{y} \right ) = w_i\left ( \underline{x} \right )y_{i,t_i}^{-1}\cdots y_{i,1}^{-1}$ and $w_{i,j}\left ( \underline{x},\underline{y} \right ) = y_{i,j}$, $j=1,\ \dots,\ t_i$) of the system $S^1$ we associate the system
$$
S_{i,j}=\left \{ E_{u_{i,j},k}\left (w_{i,j}\left ( \underline{x},\underline{y} \right ),I_A\left ( H \right )  \right )=E_{u_{i,j},k}\left (u_{i,j},I_A\left ( H \right )  \right )|\ k=1,\ \dots,\ n_{u_{i,j}}\right \}.
$$
Next, consider the system
$$
S^2=\left \{ S_{i,j}|\ i=1,\ \dots,\ m,\ j=0,\ \dots,\ t_i \right \},
$$
consisting of all the equations of all the systems $S_{i,j}.$ It is clear that if the system $S^1$ has a solution in $G$, then the system $S^2$ has a solution in $G$ too.

Since $I_A(H)$ is a non-abelian free group (by Lemma \ref{LemWSF}), $E_{u_{i,j},k}\left (w_{i,j}\left ( \underline{x},\underline{y} \right ),I_A\left ( H \right )  \right )\in I_A(F(\underline{x},\underline{y}))\cdot I_A(H)^{F(\underline{x},\underline{y})\ast I_A(H)}$ and $E_{u_{i,j},k}\left (u_{i,j},I_A\left ( H \right )  \right )\in I_A(H)$, it is clear that the system $S^2$ is of the form (\ref{Nth01}). Therefore, by Lemma \ref{SmL}, this system has a solution, say $\underline{\hat{x}}$, $\underline{\hat{y}}$, in $H$. In turn, in accordance with Condition T3) of the Main Theorem, this implies that the following equalities hold in $H/Q$
$$
\left \{ w_{i,j}\left ( \underline{\hat{x}},\underline{\hat{y}} \right ) Q=u_{i,j}Q|\ i=1,\ \dots,\ m,\ j=0,\ \dots,\ t_i\right \}.
$$
Therefore, for some elements $q_{i,j}\in Q$ we have the following set of equalities (in $H$)
$$\hat{S}^1=\left \{ w_i ( \underline{\hat{x}}  )\hat{y}_{i,t_i}^{-1}\cdots \hat{y}_{i,1}^{-1}=u_{i,0}q_{i,0},\ \hat{y}_{i,1}=u_{i,1}q_{i,1},\ \dots,\ \hat{y}_{i,t_i}=u_{i,t_i}q_{i,t_i}|\ i=1,\ \dots,\ m \right \}.$$
Whence, since $Q$ is central, there exist the elements $q_i\in Q$ such that the following equalities hold
\begin{equation}\label{SEQLm4}
\hat{S}=\left \{ w_i(\underline{\hat{x}})=u_{i,0}u_{i,1}\cdots u_{i,t_{i}}q_i=h_iq_i |\ i=1,\ \dots,\ m \right \}.
\end{equation}
This means that $\underline{\hat{x}}$ is a solution to $S$ in $H$ up to some elements $q_i\in Q$.

In accordance with (\ref{A0505SEM}), the equality $w_i(\underline{\hat{x}})=h_iq_i$ from (\ref{SEQLm4}) can be either of the form $u_i(\hat{x}_1,\dots,\hat{x}_n)=h_iq_i$ or of the form $\hat{x}_i^{m_i}u_i(\hat{x}_1,\dots,\hat{x}_n)=h_iq_i$, $m_i>0$.

Suppose that $\hat{x}_i^{m_i}u_i(\hat{x}_1,\dots,\hat{x}_n)=h_iq_i$, $m_i>0$. Then (since $Q$ is a divisible group) there is $r_i\in Q$ such that $r_i^{m_i}=q_i$.  If we replace the element $\hat{x}_i$ of the tuple $\underline{\hat{x}}$ with the element $r_i^{-1}\hat{x}_i$ (and leave the other elements unchanged), then (since $Q$ is a cental subgroup of $H$ and $u_i(x_1,\dots,x_n)\in F(x_1,\dots,x_n)'$) we have
\begin{align*}
\left ( r_i^{-1}\hat{x}_i \right )^{m_i}u_i(\hat{x}_1,\dots,r_i^{-1}\hat{x}_i,\dots,\hat{x}_n)&=r_i^{-m_i}\hat{x}_i^{m_i}u_i(\hat{x}_1,\dots,r_i^{-1}\hat{x}_i,\dots,\hat{x}_n)\\
&=q_i^{-1}\hat{x}_i^{m_i}u_i(\hat{x}_1,\dots,\hat{x}_i,\dots,\hat{x}_n)\\
&=q_i^{-1}h_iq_i=h_i,
\end{align*}
while the other equalities of $\hat{S}$ remain unchanged.

Suppose that $u_i(\hat{x}_1,\dots,\hat{x}_n)=h_iq_i$. The equation $u_i(x_1,\dots,x_n)=h_i$ has a solution in $H$ (indeed, this equation has a solution in $G$ since the system $S$ does, and $H$ is a verbally closed subgroup in $G$). It means (since $u_i(x_1,\dots,x_n)\in F(x_1,\dots,x_n)'$ ) that $h_i\in H'$, thus $q_i=h_i^{-1}u_i(\hat{x}_1,\dots,\hat{x}_n)\in H'.$ But by Condition T2) of the Main Theorem we have $Q\cap H'=\{1\}$, therefore $q_i=1$.

Thus, the tuple $(r_1^{-1}\hat{x}_1,\dots,r_l^{-1}\hat{x}_l,\hat{x}_{l+1},\dots,\hat{x}_n)$ is a solution to the system $S$ in $H$ as required.

\end{ProofMT}

\section{Acknowledgments}

The author thanks A. A. Klyachko for many useful conversations and many useful remarks.

\begin{otherlanguage}{english}

\end{otherlanguage}

\end{document}